\documentclass[a4paper,12pt,reqno]{amsart}

\usepackage[utf8]{inputenc}
\usepackage[T1]{fontenc}
\usepackage{amsmath,amssymb,amsthm,mathtools}
\usepackage{esint}
\usepackage{geometry}
\usepackage{enumitem}
\usepackage{hyperref}
\usepackage{comment}

\newtheorem{theorem}{Theorem}[section]
\newtheorem{proposition}[theorem]{Proposition}
\newtheorem{lemma}[theorem]{Lemma}

\newtheorem{example}[theorem]{Example}

\theoremstyle{definition}
\newtheorem{definition}[theorem]{Definition}

\theoremstyle{remark}
\newtheorem{remark}[theorem]{Remark}

\newcommand{\R}{\mathbb{R}}
\newcommand{\Kplus}{K_0^+}
\newcommand{\Ninf}{N_\infty}

\begin{document}
	
	\title[]{One-Sided Hardy--Littlewood Maximal Operators
		on Variable Lebesgue Spaces Beyond Log-H\"older Continuity}
	
	\author[L.~A.~Vallejos]{Lucas Alejandro Vallejos}
	\address{(L.~A.~Vallejos)  FaMAF \\ Universidad Nacional de C\'ordoba \\
		CIEM (CONICET) \\ 5000 C\'ordoba, Argentina}
	\email{lucas.vallejos@unc.edu.ar}
	
	
	\begin{abstract}
		We study the boundedness of the right-sided Hardy--Littlewood maximal
		operator on variable Lebesgue spaces over the real line. Motivated by
		the $K_0$ condition for the classical Hardy--Littlewood maximal
		operator, we introduce a directional condition $K_0^+$ adapted to the
		geometry of the one-sided maximal operator.
		
		We characterize $K_0^+$ in terms of the uniform boundedness of a
		family of one-sided averaging operators and prove that it is a
		necessary condition for the boundedness of $M^+$ on
		$L^{p(\cdot)}(\mathbb{R})$. We also establish one-sided analogues of
		the modular and reverse H\"older estimates associated with the
		classical $K_0$ condition. These estimates constitute the main local
		ingredients in the proof of boundedness of $M^+$ under the combined
		assumptions $K_0^+$ and $N_\infty$.
	\end{abstract}
	
	\date{}
	
	
	
	\maketitle
	
	\section{Introduction}
	
	Let $p(\cdot)$ be a measurable variable exponent on $\mathbb{R}$. The boundedness of the
	Hardy--Littlewood maximal operator $M$ on variable Lebesgue spaces has been extensively
	studied. Diening~\cite{Diening2004} proved that $M$ is bounded on $L^{p(\cdot)}(\Omega)$,
	$\Omega$ bounded, whenever $1<p_-\le p_+<\infty$ and $p(\cdot)$ is locally log-H\"older
	continuous ($LH_0$):
	\[
	|p(x)-p(y)| \le \frac{C_0}{\log(e+1/|x-y|)}, \qquad x,y\in\Omega .
	\]
	On the whole space, Cruz-Uribe, Fiorenza and Neugebauer~\cite{CUFN2003} showed that it
	suffices to add the log-H\"older decay condition at infinity ($LH_\infty$): there exists
	$p_\infty$ such that $|p(x)-p_\infty|\le C_\infty/\log(e+|x|)$. Independently,
	Nekvinda~\cite{Nekvinda2004} showed that $LH_\infty$ can be replaced by the weaker
	integral condition $N_\infty$ (see Definition \ref{N_infty-condition} below). These conditions are sharp
	among pointwise moduli of continuity~\cite{PickRuzicka2001,CUFN2003}. We refer to the
	monographs~\cite{CruzUribeFiorenza,DHHR-book} for a detailed account.
	
	These pointwise regularity assumptions, however, are not necessary.
	In the two-sided setting, the local log-H\"older condition can be
	replaced by the weaker $K_0$ condition. This condition admits a
	natural functional interpretation: it is equivalent to the uniform
	boundedness of the family of averaging operators over intervals or
	cubes. Combined with the $N_\infty$ condition, it yields boundedness
	of the Hardy--Littlewood maximal operator; see
	\cite{CruzUribeFiorenza}. The conditions $K_0$ and $N_\infty$ have
	also been used in the study of integral operators on variable
	Lebesgue spaces; see, for instance, \cite{UrciuoloVallejos2020}.
	
	The purpose of this paper is to develop an analogous approach for
	the one-sided Hardy--Littlewood maximal operator
	\[
	M^+f(x)
	=
	\sup_{h>0}
	\frac{1}{h}
	\int_x^{x+h}|f(y)|\,dy.
	\]
	The left-sided maximal operator is defined by
	\[
	M^-f(x)
	=
	\sup_{h>0}
	\frac{1}{h}
	\int_{x-h}^{x}|f(y)|\,dy.
	\]
	
	The theory of one-sided maximal operators has a long history in weighted
	harmonic analysis. Sawyer \cite{Sawyer1986} introduced the one-sided
	Muckenhoupt classes $A_p^+$ and $A_p^-$ and characterized the boundedness
	of the one-sided Hardy--Littlewood maximal operators on weighted Lebesgue
	spaces. Further developments of this theory were obtained by
	Mart\'in-Reyes, Ortega Salvador and de la Torre
	\cite{MartinReyesOrtegaTorre1990}, and by Mart\'in-Reyes
	\cite{MartinReyes1993}. In particular, the classes $A_\infty^+$ were
	introduced and characterized by Mart\'in-Reyes, Pick and de la Torre
	\cite{MartinReyesPickTorre1993}. One-sided reverse H\"older inequalities
	and their relationship with one-sided Muckenhoupt classes were further
	studied by Cruz-Uribe \cite{CruzUribe1997}, Cruz-Uribe, Neugebauer and
	Olesen \cite{CruzUribeNeugebauerOlesen1995}, and Riveros and de la Torre
	\cite{RiverosTorre2001}.
	
	One-sided operators have also been studied in the variable exponent
	setting. Edmunds, Kokilashvili and Meskhi
	\cite{EdmundsKokilashviliMeskhi2007,EdmundsKokilashviliMeskhi2008}
	obtained boundedness results for one-sided maximal functions, potentials
	and singular integrals on variable Lebesgue spaces under directional
	conditions on the exponent which are weaker than the usual symmetric
	log-H\"older assumptions. Weighted one- and two-weight inequalities in
	variable Lebesgue spaces were subsequently studied by Kokilashvili,
	Meskhi and Sarwar \cite{KokilashviliMeskhiSarwar2010}.
	A general self-improvement criterion for maximal operators associated
	with bases was obtained by Lerner and Ombrosi
	\cite{LernerOmbrosi2010}. Their framework includes the
	one-dimensional right-sided maximal operator. In fact, the
	one-sided basis generating $M^+$ satisfies Stein's property in
	dimension one, and this criterion will be used in the final step of
	the proof of our main theorem.
	
	The geometry of $M^+$ naturally distinguishes intervals lying to
	the left from intervals lying to the right: for $a<b<c$, the averages
	of $|f|$ over $(b,c)$ are controlled by $M^+f$ at every point of
	$(a,b)$. This leads to a directional version of the $K_0$ condition,
	which we denote by $K_0^+$; see Definition~\ref{Kplus}.
	
	Our main result is the following.
	
	\begin{theorem}\label{main-Theorem}
		Let $p(\cdot)\in\mathcal{P}(\mathbb{R})$ satisfy $1<p_-\leq p_+<\infty.$ Assume that $p(\cdot)\in K_0^+(\mathbb{R})\cap N_\infty(\mathbb{R}).$
		Then the right-sided Hardy--Littlewood maximal operator is bounded on
		$L^{p(\cdot)}(\mathbb{R})$; that is,
		\[
		\|M^+f\|_{p(\cdot)}
		\leq
		C\|f\|_{p(\cdot)}
		\]
		for every $f\in L^{p(\cdot)}(\mathbb{R})$.
	\end{theorem}
	
	Since $LH_0\cap LH_\infty\subset K_0\cap N_\infty$ (see~\cite[Propositions~4.9 and~4.57]{CruzUribeFiorenza}) and
	$K_0\subset K_0^+$ by~\eqref{examples-K0-inclusion}, Theorem~\ref{main-Theorem}
	contains the log-H\"older case. We note that in this case the boundedness of $M^+$
	already follows from that of $M$, since $M^+f\le Mf$ pointwise.
	Example~\ref{strict-one-sided} shows that Theorem~\ref{main-Theorem} applies to
	exponents that are not covered by the two-sided theory, even under $N_\infty$.
	
	The main new ingredients in the proof are one-sided counterparts of
	the modular and reverse H\"older estimates associated with the
	classical $K_0$ condition; see Lemmas~\ref{lemma455-CU}
	and~\ref{lemma456-CU}.
	
	The paper is organized as follows. In Section~2, we collect
	the known preliminary results used throughout the paper.
	In Section~3, we introduce the condition $K_0^+$ and establish
	equivalent formulations in terms of certain averaging operators.
	In Section~4, we prove that $K_0^+$ is necessary for the
	boundedness of $M^+$ on $L^{p(\cdot)}(\mathbb R)$.
	Section~5 is devoted to the proof of our main result,
	Theorem~\ref{main-Theorem}, preceded by four auxiliary lemmas
	that provide the key ingredients of the argument.
	Finally, in Section~6, we present examples of exponents
	satisfying the condition $K_0^+$.
	
	\section{Preliminaries}
	
	Let $\mathcal{P}(\mathbb{R})$ denote the set of measurable functions $p(\cdot):\mathbb{R}\longrightarrow[1,\infty).$
	For a measurable set $\Omega\subset\mathbb R$ of positive measure, we write
	\[
	p_-(\Omega)=\operatorname*{ess\,inf}_{x\in\Omega}p(x),
	\qquad
	p_+(\Omega)=\operatorname*{ess\,sup}_{x\in\Omega}p(x).
	\]
	When $\Omega=\mathbb R$ we simply write $p_-=p_-(\mathbb R)$ and $p_+=p_+(\mathbb R)$.
	
	Throughout the paper we assume $1<p_-\leq p_+<\infty.$
	The conjugate exponent $p'(\cdot)$ is defined pointwise by
	\[
	\frac{1}{p(x)}+\frac{1}{p'(x)}=1.
	\]
	
	The modular associated with $p(\cdot)$ is
	\[
	\rho_{p(\cdot)}(f)
	=
	\int_{\mathbb{R}}
	|f(x)|^{p(x)}\,dx.
	\]
	
	Given $p(\cdot)\in \mathcal{P}(\mathbb{R})$, let $L^{p(\cdot)}(\mathbb{R})$ be the Banach space of measurable functions $f$ on $\mathbb{R}$ such that, for some  $\lambda >0$,  the modular $\rho_{p(\cdot)}(f/\lambda)<\infty$ is finite, with norm
	\begin{equation*}
		\left\Vert f\right\Vert _{p(\cdot )}=\inf \{\lambda >0:~%
		\mathit{\rho }_{p(\cdot )}\left( f/\lambda \right) \leq 1\}\text{.}
	\end{equation*}
	
	These spaces are known as {\it variable Lebesgue spaces} and generalize the classical Lebesgue spaces $L^p(\R)$.
	
	We recall several basic properties of variable Lebesgue spaces. See, for example: \cite{CruzUribeFiorenza}.
	
	\begin{lemma}
		\label{norm-modular-boundedness}
		Let $p(\cdot):\mathbb R\to[1,\infty)$ be a measurable
		exponent with $p_+<\infty$.
		A bound on $\|f\|_{p(\cdot)}$ yields a bound on
		$\rho_{p(\cdot)}(f)$, and conversely.
		More precisely, for $C_1,C_2>0$,
		\[
		\|f\|_{p(\cdot)}\leq C_1
		\quad\Longrightarrow\quad
		\rho_{p(\cdot)}(f)
		\leq \max\{C_1^{p_-},C_1^{p_+}\},
		\]
		and
		\[
		\rho_{p(\cdot)}(f)\leq C_2
		\quad\Longrightarrow\quad
		\|f\|_{p(\cdot)}
		\leq \max\{C_2^{1/p_-},C_2^{1/p_+}\}.
		\]
		In particular,
		\[
		\|f\|_{p(\cdot)}\leq1
		\quad\Longleftrightarrow\quad
		\rho_{p(\cdot)}(f)\leq1.
		\]
	\end{lemma}
	
	\begin{lemma}
		\label{comparision-modular-norm}
		Let $p(\cdot):\mathbb R\to[1,\infty)$ be a measurable
		exponent with $p_+<\infty$, and let
		$\Omega\subset\mathbb R$ be measurable with
		$|\Omega|>0$.
		Then the following assertions hold.
		\begin{enumerate}[label=\textnormal{(\alph*)}]
			\item If $\|f\chi_\Omega\|_{p(\cdot)}\leq1$, then
			\[
			\|f\chi_\Omega\|_{p(\cdot)}^{p_+(\Omega)}
			\leq
			\int_\Omega |f(x)|^{p(x)}\,dx
			\leq
			\|f\chi_\Omega\|_{p(\cdot)}^{p_-(\Omega)}.
			\]
			
			\item If $\|f\chi_\Omega\|_{p(\cdot)}\geq1$, then
			\[
			\|f\chi_\Omega\|_{p(\cdot)}^{p_-(\Omega)}
			\leq
			\int_\Omega |f(x)|^{p(x)}\,dx
			\leq
			\|f\chi_\Omega\|_{p(\cdot)}^{p_+(\Omega)}.
			\]
		\end{enumerate}
		In particular, if $\|f\|_{p(\cdot)}\leq1$, then
		\[
		\int_{\mathbb R}|f(x)|^{p(x)}\,dx
		\leq \|f\|_{p(\cdot)}.
		\]
	\end{lemma}
	
	\begin{lemma}
		\label{unit-norm-modular}
		Let $p(\cdot):\mathbb R\to[1,\infty)$ be a measurable
		exponent with $p_+<\infty$.
		If $\|f\|_{p(\cdot)}=1$, then
		\[
		\int_{\mathbb R}|f(x)|^{p(x)}\,dx=1.
		\]
	\end{lemma}
	
	As a consequence of Lemma~\ref{unit-norm-modular},
	for every nonzero $f\in L^{p(\cdot)}(\mathbb R)$,
	\[
	\rho_{p(\cdot)}
	\left(\frac{f}{\|f\|_{p(\cdot)}}\right)=1.
	\]
	In particular, for every measurable set
	$E\subset\mathbb R$ with $0<|E|<\infty$,
	\begin{equation}\label{char-normalization}
		\int_E
		\left(\frac{1}{\|\chi_E\|_{p(\cdot)}}\right)^{p(x)}
		\,dx=1.
	\end{equation}
	
	
	\begin{lemma}\label{dilationExpo}
		Let $p(\cdot )\in \mathcal{P}(\mathbb{R})$. Then for all $s$ with $1/p_{-}\leq s<\infty $,
		\begin{equation*}
			\Vert |f|^{s}\Vert _{p(\cdot)}=\Vert f\Vert_{sp(\cdot)}^{s}.
		\end{equation*}
	\end{lemma}
	
	\begin{lemma}\label{Holder-variable}
		Let $p(\cdot)$ be a measurable exponent and let
		$p'(\cdot)$ be its pointwise conjugate.
		There exists a constant $C>0$ such that
		\begin{equation*}
			\int_{\mathbb R}|f(x)g(x)|\,dx \leq
			C\|f\|_{p(\cdot)}\|g\|_{p'(\cdot)}
		\end{equation*}
		for every $f\in L^{p(\cdot)}(\mathbb R)$ and
		$g\in L^{p'(\cdot)}(\mathbb R)$.
	\end{lemma}
	
	
	
	\begin{lemma}
		\label{associate-norm}
		Let $p(\cdot)$ be a measurable exponent.
		There exist constants $c,C>0$, independent of $f$,
		such that for every $f\in L^{p(\cdot)}(\mathbb R)$
		there exists $g\in L^{p'(\cdot)}(\mathbb R)$ with
		$\|g\|_{p'(\cdot)}\leq1$ satisfying
		\[
		c\|f\|_{p(\cdot)}
		\leq
		\left|\int_{\mathbb R}f(x)g(x)\,dx\right|
		\leq
		C\|f\|_{p(\cdot)}.
		\]
	\end{lemma}

	\subsection*{The condition at infinity}
	
	We retain the usual $N_\infty$ condition.
	
	\begin{definition}\label{N_infty-condition}
		We say that
		\[
		p(\cdot)\in N_\infty(\mathbb{R})
		\]
		if there exist constants $\Lambda_\infty>0$ and
		$p_\infty\in[1,\infty)$ such that
		\[
		\int_{\Omega_\infty}
		\exp
		\left(
		-\Lambda_\infty
		\left|
		\frac{1}{p(x)}
		-
		\frac{1}{p_\infty}
		\right|^{-1}
		\right)
		dx
		<\infty,
		\]
		where
		\[
		\Omega_\infty
		=
		\left\{
		x\in\mathbb{R}:
		\left|
		\frac{1}{p(x)}
		-
		\frac{1}{p_\infty}
		\right|>0
		\right\}.
		\]
	\end{definition}
	
	Since $p_+<\infty$, this condition is equivalent, after changing the
	constant, to
	\[
	\int_{\Omega_\infty}
	\exp
	\left(
	-\Lambda_\infty
	|p(x)-p_\infty|^{-1}
	\right)
	dx
	<\infty.
	\]
	
	\subsection*{One-sided weights}
	
	A weight is a non-negative locally integrable function
	$w:\mathbb{R}\to[0,\infty)$. For a measurable set $E$, we write
	\[
	w(E)=\int_E w(x)\,dx.
	\]
	
	We shall use a classical one-sided $A_\infty$ principle. A convenient
	formulation is the following.
	
	\begin{definition}\label{Ainfplus-property}
		We say that a weight $w$ satisfies the $A_\infty^+$ property if,
		for every $\eta\in(0,1)$, there exists
		$\beta=\beta(\eta)\in(0,1)$ such that for all $a<b<c$ and every measurable set $E\subset(b,c)$, $w(E)\leq\beta w(a,b)$ implies
		$\frac{|E|}{c-a}<\eta.$
		
	\end{definition}
	
	This is one of the standard equivalent formulations of the
	one-sided $A_\infty$ condition; see
	\cite{MartinReyesPickTorre1993,CruzUribe1997}.
	
	We shall also use the one-sided reverse H\"older classes.
	
	\begin{definition}
		Let $r>1$. We say that a weight $w$ satisfies the one-sided reverse
		H\"older condition $RH_r^+$ if there exists $C>0$ such that for
		every pair of adjacent intervals of equal length,
		\[
		I^-=(a-h,a),
		\qquad
		I^+=(a,a+h),
		\]
		one has
		\begin{equation}\label{RHplus}
			\left(
			\frac{1}{|I^-|}
			\int_{I^-}w(x)^r\,dx
			\right)^{1/r}
			\leq
			C
			\frac{1}{|I^+|}
			\int_{I^+}w(x)\,dx.
		\end{equation}
	\end{definition}
	
	The formulation \eqref{RHplus} is equivalent to the usual
	definition of $RH_r^+$; see \cite{RiverosTorre2001}.
	
	\subsection*{A level-set estimate for $M^+$}
	
	We shall use the following estimate due to Mart\'in-Reyes,
	Rivera-R\'ios and Rodr\'iguez-Padilla.
	
	\begin{lemma}[{\cite[Lemma~4]{MartinReyesRiveraRiosRodriguezPadilla2025}}]
		\label{sparse-alike}
		Let $v\in L^1_{\mathrm{loc}}(\mathbb R)$ and $0<\lambda_1<\lambda_2$, and set
		\[
		F=\{z\in\mathbb R:M^+v(z)\leq\lambda_2\}.
		\]
		If $x\in\mathbb R$ satisfies $M^+v(x)\leq\lambda_1$, then for every $y>x$,
		\[
		|F\cap(x,y)|
		\geq
		\left(1-\frac{\lambda_1}{\lambda_2}\right)(y-x).
		\]
	\end{lemma}

	\section{The condition \texorpdfstring{$K_0^+$}{K0+}}
	
	Our definition is also motivated by the classical $K_0$ condition in
	variable Lebesgue spaces. The $K_0$ condition was first considered in
	the more general framework of Banach function spaces by Berezhnoi
	\cite{Berezhnoi1999}, and its importance in variable Lebesgue spaces
	was emphasized by Kopaliani \cite{Kopaliani2007}. The connection between
	$K_0$, averaging operators, modular estimates and the Hardy--Littlewood
	maximal operator is developed in detail by Cruz-Uribe and Fiorenza
	\cite{CruzUribeFiorenza}; see also Diening \cite{Diening2005} and
	Lerner \cite{Lerner2010}.
	
	The directional nature of $K_0^+$ is also consistent with other
	non-symmetric Muckenhoupt-type conditions arising in harmonic analysis.
	In particular, classes involving weights evaluated on geometrically
	different sets or through invertible transformations appear in the work
	of Iba\~nez-Firnkorn, Riveros and Vidal \cite{IbanezRiverosVidal2022};
	see also the survey by Riveros and Vidal \cite{RiverosVidal2023}.

	\begin{definition}\label{Kplus}
		Let $p(\cdot)\in\mathcal{P}(\mathbb{R})$. We say that $p(\cdot)\in K_0^+(\mathbb{R})$ if
		\begin{equation*}
			[p]_{K_0^+}
			:=
			\sup_{a<b<c}
			\frac{
				\|\chi_{(a,b)}\|_{p(\cdot)}
				\|\chi_{(b,c)}\|_{p'(\cdot)}
			}{
				c-a
			}
			<\infty.
		\end{equation*}
	\end{definition}
	
	The orientation in Definition \ref{Kplus} is essential. The
	$L^{p(\cdot)}$ norm is taken on the interval lying to the left and
	the $L^{p'(\cdot)}$ norm on the interval lying to the right.
	
	The following observation shows that Definition
	\ref{Kplus} is equivalent to the apparently simpler definition
	using adjacent intervals of equal length.
	
	\begin{proposition}\label{Kplus-equivalent}
		Let $p(\cdot)\in\mathcal{P}(\mathbb{R})$. Then $p(\cdot)\in K_0^+(\mathbb{R})$ if and only if
		\begin{equation}\label{Kplus-equal}
			\sup_{\substack{a\in\mathbb{R}\\h>0}}
			\frac{
				\|\chi_{(a-h,a)}\|_{p(\cdot)}
				\|\chi_{(a,a+h)}\|_{p'(\cdot)}
			}{h}
			<\infty.
		\end{equation}
	\end{proposition}
	
	\begin{proof}
		It is immediate that $p(\cdot)\in K_0^+(\mathbb{R})$ implies
		\eqref{Kplus-equal}, up to an absolute factor, by taking $a-h<a<a+h.$
		
		Conversely, suppose that \eqref{Kplus-equal} holds and let
		$a<b<c$. Set $\ell_-=b-a, \ \ell_+=c-b.$
		
		If $\ell_-\leq\ell_+$, then $(a,b)\subset(b-\ell_+,b),$ and hence
		\[
		\|\chi_{(a,b)}\|_{p(\cdot)}
		\leq
		\|\chi_{(b-\ell_+,b)}\|_{p(\cdot)}.
		\]
		Therefore,
		\begin{equation*}
			\|\chi_{(a,b)}\|_{p(\cdot)}
			\|\chi_{(b,c)}\|_{p'(\cdot)} \leq
			C\ell_+\leq C(c-a).
		\end{equation*}
		
		If $\ell_+\leq\ell_-$, we enlarge the right interval instead: $(b,c)\subset(b,b+\ell_-),$ and obtain
		\[
		\|\chi_{(a,b)}\|_{p(\cdot)}
		\|\chi_{(b,c)}\|_{p'(\cdot)}
		\leq
		C\ell_-
		\leq
		C(c-a).
		\]
		This proves that $p(\cdot) \in K_0^+(\mathbb{R})$.
	\end{proof}
	
	\subsection*{One-sided averaging operators}
	
	Given $a<b<c$, we define
	\begin{equation}\label{Aabc}
		A_{a,b,c}^+f(x)
		=
		\frac{\chi_{(a,b)}(x)}{c-a}
		\int_b^c f(y)\,dy.
	\end{equation}
	
	These operators reflect exactly the geometry of $M^+$.
	
	\begin{proposition}\label{averaging-characterization}
		The following statements are equivalent.
		\begin{enumerate}[label=\textnormal{(\roman*)}]
			\item
			$p(\cdot)\in K_0^+(\mathbb{R})$;
			
			\item
			the family $\left\{A_{a,b,c}^+:a<b<c\right\}$ is uniformly bounded on $L^{p(\cdot)}(\mathbb{R})$.
		\end{enumerate}
	\end{proposition}
	
	\begin{proof}
		Suppose first that
		$p(\cdot)\in K_0^+(\mathbb{R})$.
		By H\"older's inequality, Lemma \ref{Holder-variable},
		\begin{align*}
			\|A_{a,b,c}^+f\|_{p(\cdot)}
			&=
			\frac{1}{c-a}
			\left|
			\int_b^c f(y)\,dy
			\right|
			\|\chi_{(a,b)}\|_{p(\cdot)}
			\\
			&\leq
			C_{p(\cdot)}
			\frac{
				\|\chi_{(a,b)}\|_{p(\cdot)}
				\|\chi_{(b,c)}\|_{p'(\cdot)}
			}{
				c-a
			}
			\|f\chi_{(b,c)}\|_{p(\cdot)}
			\\
			&\leq
			C[p]_{K_0^+}
			\|f\|_{p(\cdot)}.
		\end{align*}
		Thus the family is uniformly bounded.
		
		Conversely, suppose that
		\[
		\|A_{a,b,c}^+f\|_{p(\cdot)}
		\leq
		C\|f\|_{p(\cdot)}
		\]
		with $C$ independent of $a,b,c$. By the associate norm, Lemma \ref{associate-norm}, of $L^{p'(\cdot)}$,
		\[
		\|\chi_{(b,c)}\|_{p'(\cdot)}
		\leq
		C
		\sup_{\|f\|_{p(\cdot)}\leq1}
		\int_b^c |f(y)|\,dy.
		\]
		Hence
		\begin{align*}
			&
			\|\chi_{(a,b)}\|_{p(\cdot)}
			\|\chi_{(b,c)}\|_{p'(\cdot)}
			\\
			&\qquad\leq
			C
			\sup_{\|f\|_{p(\cdot)}\leq1}
			\|\chi_{(a,b)}\|_{p(\cdot)}
			\int_b^c |f(y)|\,dy
			\\
			&\qquad=
			C(c-a)
			\sup_{\|f\|_{p(\cdot)}\leq1}
			\|A_{a,b,c}^+|f|\|_{p(\cdot)}
			\\
			&\qquad\leq
			C(c-a).
		\end{align*}
		Therefore $p(\cdot)\in K_0^+(\mathbb{R}).$
		
	\end{proof}
	
	\section{Necessity of the \texorpdfstring{$K_0^+$}{K0+} condition}
	
	We now show that $K_0^+$ is a necessary condition for the
	boundedness of $M^+$.
	
	\begin{theorem}\label{necessity}
		Suppose that $M^+:L^{p(\cdot)}(\mathbb{R}) \longrightarrow
		L^{p(\cdot)}(\mathbb{R})$ is bounded. Then $p(\cdot)\in K_0^+(\mathbb{R}).$
		
	\end{theorem}
	
	\begin{proof}
		Fix $a<b<c$. For every $x\in(a,b)$, the interval $[x,c]$
		contains $(b,c)$ and $c-x\leq c-a.$
		Consequently,
		\[
		M^+f(x)
		\geq
		\frac{1}{c-a}
		\int_b^c |f(y)|\,dy.
		\]
		Thus
		\[
		A_{a,b,c}^+|f|(x)
		\leq
		M^+f(x),
		\qquad
		x\in(a,b).
		\]
		
		Therefore,
		\[
		\|A_{a,b,c}^+|f|\|_{p(\cdot)}
		\leq
		\|M^+f\|_{p(\cdot)}
		\leq
		C\|f\|_{p(\cdot)}.
		\]
		The constant is independent of $a,b,c$. Hence the family of
		one-sided averaging operators is uniformly bounded. Proposition
		\ref{averaging-characterization} now yields $p(\cdot)\in K_0^+(\mathbb{R}).$
		
	\end{proof}
	
	\begin{remark}
		Theorem \ref{necessity} provides a functional justification for
		Definition \ref{Kplus}. The condition $K_0^+$ is not introduced
		ad hoc: it is precisely the characteristic-function condition
		associated with the one-sided averaging operators naturally
		dominated by $M^+$.
	\end{remark}
	
	\section{Proof of the main result}
	
	We now develop the two local estimates that replace the corresponding
	$K_0$ lemmas in the two-sided theory.
	
	The following result is the one-sided counterpart of
	\cite[Lemma~4.55]{CruzUribeFiorenza}.
	In the classical setting, that lemma shows that the $K_0$ condition
	implies a local modular inequality for averaging operators.
	Here the averaging and testing intervals are separated according to the
	geometry of $M^+$. The proof follows the general strategy of
	\cite[Lemma~4.55]{CruzUribeFiorenza}, with the essential modification
	that the usual $K_0$ estimate is replaced by the directional
	$K_0^+$ condition.
	
	\begin{lemma}\label{lemma455-CU}
		Let $p(\cdot)\in\mathcal{P}(\mathbb{R})$ satisfy $1<p_-\leq p_+<\infty$
		and suppose $p(\cdot)\in K_0^+(\mathbb{R})$.
		Let $a<b<c$ and let
		$f\in L^{p(\cdot)}(\mathbb{R})$. Assume that there exist constants
		$c_1,c_2>0$ such that $\frac{1}{c-a} \int_b^c |f(y)|\,dy \geq c_1$ and $\|f\|_{p(\cdot)} \leq c_2.$ Then
		\begin{equation}\label{455plus}
			\int_a^b
			\left(
			\frac{1}{c-a}
			\int_b^c |f(y)|\,dy
			\right)^{p(x)}
			dx
			\leq
			C
			\int_b^c |f(y)|^{p(y)}\,dy,
		\end{equation}
		where $C$ depends only on
		$p(\cdot)$, $c_1$, $c_2$, and $[p]_{K_0^+}$.
	\end{lemma}
	
	\begin{proof}
		We adapt the argument of
		\cite[Lemma~4.55]{CruzUribeFiorenza} to the one-sided geometry.
		The main difference is that the variable $x$ belongs to the left
		interval $(a,b)$, whereas the averaging variable $y$ belongs to the right interval $(b,c)$. Set $D=c-a.$
		Since $p'_-=(p_+)'>1,$ the function
		\[
		\Phi(\alpha)
		=
		\frac{1}{D}
		\int_b^c \alpha^{p'(y)-1}\,dy,
		\qquad \alpha>0,
		\]
		is continuous and strictly increasing. Moreover,
		\[
		\Phi(\alpha)\longrightarrow0
		\quad\text{as }\alpha\to0^+,
		\qquad
		\Phi(\alpha)\longrightarrow\infty
		\quad\text{as }\alpha\to\infty.
		\]
		Hence there exists a unique $\alpha>0$ such that
		\begin{equation}\label{alpha-definition}
			\frac{1}{D}
			\int_b^c
			\alpha^{p'(y)-1}\,dy
			=
			\frac{1}{D}
			\int_b^c |f(y)|\,dy.
		\end{equation}
		
		Set
		\[
		B=
		\frac{1}{D}
		\int_b^c \alpha^{p'(y)}\,dy.
		\]
		Since $p'(x)(p(x)-1)=p(x),$ we obtain
		
		\begin{equation}\label{455-factorization}
			\int_a^b
			\left(
			\frac{1}{D}
			\int_b^c |f(y)|\,dy
			\right)^{p(x)}
			dx
			=
			B
			\int_a^b
			\left(
			\frac{1}{D}
			\int_b^c
			\alpha^{p'(y)-p'(x)}
			\,dy
			\right)^{p(x)-1}
			dx.
		\end{equation}
		
		We estimate both factors. By \eqref{alpha-definition},
		\[
		\int_b^c\alpha^{p'(y)}\,dy
		=
		\alpha
		\int_b^c|f(y)|\,dy.
		\]
		Therefore,
		\[
		\int_b^c\alpha^{p'(y)}\,dy
		=
		2\alpha\int_b^c|f(y)|\,dy
		-
		\int_b^c\alpha^{p'(y)}\,dy.
		\]
		It follows that
		\begin{align*}
			\int_b^c\alpha^{p'(y)}\,dy
			&\leq
			2\alpha
			\int_{\{
				y\in(b,c):
				2\alpha|f(y)|>\alpha^{p'(y)}
				\}}
			|f(y)|\,dy.
		\end{align*}
		
		On the set appearing on the right-hand side, $2|f(y)|>\alpha^{p'(y)-1},$ and hence 
		
		$\alpha<(2|f(y)|)^{p(y)-1}.$ Consequently,
		\[
		2\alpha|f(y)|
		\leq
		2^{p(y)}|f(y)|^{p(y)}
		\leq
		2^{p_+}|f(y)|^{p(y)}.
		\]
		Thus
		\begin{equation}\label{B-estimate}
			B
			\leq
			2^{p_+}
			\frac{1}{D}
			\int_b^c
			|f(y)|^{p(y)}\,dy.
		\end{equation}
		
		We now estimate the second factor in
		\eqref{455-factorization}. By our hypotheses and
		\eqref{alpha-definition},
		\[
		c_1
		\leq
		\frac{1}{D}
		\int_b^c
		\alpha^{p'(y)-1}\,dy.
		\]
		Since $c-b\leq D$, if $0<\alpha\leq1$, then
		\[
		c_1
		\leq
		\frac1D\int_b^c\alpha^{p'(y)-1}\,dy
		\leq
		\alpha^{p'_- -1}.
		\]
		Thus, in all cases,
		\begin{equation}\label{alpha-lower}
			\alpha\geq
			d_1:=
			\min\left\{
			1,
			c_1^{1/(p'_--1)}
			\right\}>0.
		\end{equation}
		
		On the other hand, by the inequality H\"older and our hypothesis for $f$,
		\begin{align*}
			\int_b^c |f(y)|\,dy
			&\leq
			C_{p(\cdot)}
			\|f\|_{p(\cdot)}
			\|\chi_{(b,c)}\|_{p'(\cdot)}
			\\
			&\leq
			c_2C_{p(\cdot)}
			\|\chi_{(b,c)}\|_{p'(\cdot)}.
		\end{align*}
		
		Let $N=\|\chi_{(b,c)}\|_{p'(\cdot)}.$
		Applying \eqref{char-normalization} to the exponent $p'(\cdot)$,
		we have
		\[
		\int_b^c
		N^{-p'(y)}\,dy
		=
		1.
		\]
		Choose $d_2\geq1$ sufficiently large so that
		\[
		d_2^{p'_--1}
		\geq
		c_2C_{p(\cdot)}.
		\]
		Then
		\begin{align*}
			\int_b^c
			\left(
			\frac{d_2}{N}
			\right)^{p'(y)-1}
			dy
			&=
			N
			\int_b^c
			d_2^{p'(y)-1}
			N^{-p'(y)}
			dy
			\\
			&\geq
			c_2C_{p(\cdot)}N
			\\
			&\geq
			\int_b^c|f(y)|\,dy.
		\end{align*}
		By \eqref{alpha-definition} and monotonicity,
		\begin{equation}\label{alpha-upper}
			\alpha
			\leq
			\frac{d_2}{N}.
		\end{equation}
		
		For each $x\in(a,b)$, define
		\[
		E_+(x)
		=
		\{y\in(b,c):p'(y)>p'(x)\}
		\]
		and
		\[
		E_-(x)
		=
		\{y\in(b,c):p'(y)\leq p'(x)\}.
		\]
		
		Using \eqref{alpha-lower} on $E_-(x)$ and
		\eqref{alpha-upper} on $E_+(x)$, we obtain
		\begin{align*}
			&
			\int_b^c
			\alpha^{p'(y)-p'(x)}\,dy
			\\
			&\leq
			D_1
			N^{p'(x)}
			\int_b^c
			N^{-p'(y)}\,dy
			+
			D_2(c-b)
			\\
			&=
			D_1N^{p'(x)}
			+
			D_2(c-b),
		\end{align*}
		where $D_1,D_2$ depend only on
		$p(\cdot)$, $d_1$, and $d_2$.
		
		Hence
		\[
		\frac{1}{D}
		\int_b^c
		\alpha^{p'(y)-p'(x)}\,dy
		\leq
		D_1
		\frac{N^{p'(x)}}{D}
		+
		D_2.
		\]
		
		Therefore,
		\begin{align}
			&
			\int_a^b
			\left(
			\frac{1}{D}
			\int_b^c
			\alpha^{p'(y)-p'(x)}
			\,dy
			\right)^{p(x)-1}
			dx
			\nonumber
			\\
			&\leq
			C D
			\int_a^b
			\left(
			\frac{N}{D}
			\right)^{p(x)}
			dx
			+
			C(b-a).
			\label{first-factor-before-K}
		\end{align}
		
		Since $p(\cdot)\in K_0^+(\mathbb{R})$,
		\[
		\frac{
			\|\chi_{(a,b)}\|_{p(\cdot)}
			N
		}{D}
		\leq
		[p]_{K_0^+}.
		\]
		Thus
		\[
		\frac{N}{D}
		\leq
		\frac{
			[p]_{K_0^+}
		}{
			\|\chi_{(a,b)}\|_{p(\cdot)}
		}.
		\]
		Using \eqref{char-normalization},
		\begin{align*}
			D
			\int_a^b
			\left(
			\frac{N}{D}
			\right)^{p(x)}
			dx
			&\leq
			CD
			\int_a^b
			\left(
			\frac{1}{
				\|\chi_{(a,b)}\|_{p(\cdot)}
			}
			\right)^{p(x)}
			dx
			\\
			&=
			CD.
		\end{align*}
		Since $b-a\leq D$, \eqref{first-factor-before-K}
		gives
		\begin{equation}\label{first-factor-final}
			\int_a^b
			\left(
			\frac{1}{D}
			\int_b^c
			\alpha^{p'(y)-p'(x)}
			\,dy
			\right)^{p(x)-1}
			dx
			\leq
			CD.
		\end{equation}
		
		Combining
		\eqref{455-factorization},
		\eqref{B-estimate}, and
		\eqref{first-factor-final}, we obtain
		\[
		\int_a^b
		\left(
		\frac{1}{c-a}
		\int_b^c |f(y)|\,dy
		\right)^{p(x)}
		dx
		\leq
		C
		\int_b^c |f(y)|^{p(y)}\,dy.
		\]
		This proves the lemma.
	\end{proof}
	
	The following lemma follows the argument in
	\cite[Theorem~1, $(c)\Rightarrow(d)\Rightarrow(e)\Rightarrow(f)$]
	{MartinReyesPickTorre1993}, with $g\equiv1$.
	We give the proof to make explicit that the condition is
	required only for triples contained in a fixed interval.
	
	\begin{lemma}\label{local-one-sided-RH}
		Let $U=(u,v)$ be a bounded interval, and let
		$w\in L^\infty(U)$ satisfy $w>0$ almost everywhere in $U$.
		Suppose that, for every $\eta\in(0,1)$, there exists
		$\beta_\eta\in(0,1)$ such that
		\begin{equation}\label{local-weight-condition}
			w(E)\leq\beta_\eta w(a,b)
			\quad\Longrightarrow\quad
			|E|<\eta(c-a)
		\end{equation}
		for every $u\leq a<b<c\leq v$ and every measurable
		set $E\subset(b,c)$.
		
		Then there exist $r>1$ and $C>0$ such that
		\begin{equation}\label{local-weight-RH}
			\left(
			\fint_{a-h}^{a}w(x)^r\,dx
			\right)^{1/r}
			\leq
			C\fint_a^{a+h}w(x)\,dx
		\end{equation}
		whenever $h>0$ and
		\[
		u\leq a-h<a<a+h\leq v.
		\]
		The constants $r$ and $C$ depend only on
		$\beta_{1/4}$ and $\beta_{1/2}$.
		In particular, they are independent of $U$ and
		$\|w\|_{L^\infty(U)}$.
	\end{lemma}
	
	\begin{proof}
		We follow the proof of
		\cite[Theorem~1, $(c)\Rightarrow(d)\Rightarrow(e)\Rightarrow(f)$]
		{MartinReyesPickTorre1993}, with $g\equiv1$, keeping all
		intervals inside $U$.
		
		Set
		\[
		\kappa=\frac{\beta_{1/4}}4.
		\]
		We first show that
		\begin{equation}\label{weight-density-local}
			|\{x\in J:w(x)>\kappa\lambda\}|
			\geq\frac12|J|
		\end{equation}
		whenever $J=(a,b)\subset U$ and $\lambda>0$ satisfy
		\[
		\frac{w(a,x)}{x-a}\geq\lambda
		\qquad(a<x<b).
		\]
		
		Let $x_0=b$. For each negative integer $k$, choose
		$x_k\in(a,x_{k+1})$ so that
		\[
		w(a,x_k)=\frac12w(a,x_{k+1}).
		\]
		These points exist because $w>0$ almost everywhere
		and its integral is continuous.
		Moreover,
		\[
		w(a,x_k)=2^k w(a,b),
		\]
		and hence $x_k\to a$ as $k\to-\infty$.
		For $k\leq-1$, define
		\[
		I_k=(x_{k-1},x_k),
		\qquad
		E=\{x\in J:w(x)\leq\kappa\lambda\},
		\qquad
		E_k=E\cap(x_k,x_{k+1}).
		\]
		The construction gives
		\[
		w(a,x_{k+1})=4w(I_k).
		\]
		By the hypothesis on the averages, also at $x=b$
		by passage to the limit,
		\[
		\lambda(x_{k+1}-a)
		\leq w(a,x_{k+1}).
		\]
		Consequently,
		\[
		w(E_k)
		\leq\kappa\lambda|E_k|
		\leq\kappa\lambda(x_{k+1}-a)
		\leq4\kappa w(I_k)
		=\beta_{1/4}w(I_k).
		\]
		Applying \eqref{local-weight-condition} to
		$x_{k-1}<x_k<x_{k+1}$ yields
		\[
		|E_k|<\frac14(x_{k+1}-x_{k-1}).
		\]
		All these triples lie in $[a,b]\subset[u,v]$.
		Since the intervals $(x_k,x_{k+1})$, $k\leq-1$,
		partition $J$ up to a countable set, we obtain
		\[
		|E|
		=\sum_{k=-\infty}^{-1}|E_k|
		\leq\frac14
		\sum_{k=-\infty}^{-1}(x_{k+1}-x_{k-1})
		\leq\frac12|J|.
		\]
		This proves \eqref{weight-density-local}.
		
		Fix now $J=(a,b)\subset U$, extend $w\chi_J$ by zero
		outside $J$, and set
		\[
		\lambda_0=M^-(w\chi_J)(b).
		\]
		For $\lambda>\lambda_0$, write
		\[
		\{x\in\mathbb R:M^-(w\chi_J)(x)>\lambda\}
		=
		\bigcup_j J_j,
		\]
		where $J_j=(a_j,b_j)$ are its pairwise disjoint
		component intervals.
		As in the proof of $(e)\Rightarrow(f)$ in the cited
		theorem, these intervals are contained in $J$ and satisfy
		\[
		\frac{w(a_j,b_j)}{b_j-a_j}=\lambda,
		\qquad
		\frac{w(a_j,x)}{x-a_j}\geq\lambda
		\quad(a_j<x<b_j).
		\]
		Thus \eqref{weight-density-local} applies to every $J_j$.
		By the Lebesgue differentiation theorem,
		\begin{align*}
			w(\{x\in J:w(x)>\lambda\})
			&\leq\sum_j w(J_j)\\
			&=\lambda\sum_j|J_j|\\
			&\leq
			2\lambda
			\sum_j|\{x\in J_j:w(x)>\kappa\lambda\}|\\
			&\leq
			2\lambda
			|\{x\in J:w(x)>\kappa\lambda\}|.
		\end{align*}
		
		Multiply this inequality by $\lambda^{\delta-1}$,
		where $\delta>0$, and integrate over
		$(\lambda_0,\infty)$.
		Tonelli's theorem gives
		\[
		\frac1\delta
		\int_{\{x\in J:w(x)>\lambda_0\}}
		w(x)\bigl(w(x)^\delta-\lambda_0^\delta\bigr)\,dx
		\leq
		\frac{2}
		{(1+\delta)\kappa^{1+\delta}}
		\int_J w(x)^{1+\delta}\,dx.
		\]
		The integrand
		$w(x)(w(x)^\delta-\lambda_0^\delta)$ is non-positive
		where $w(x)\leq\lambda_0$.
		Therefore,
		\[
		\frac1\delta
		\left(
		\int_Jw(x)^{1+\delta}\,dx
		-\lambda_0^\delta w(J)
		\right)
		\leq
		\frac{2}
		{(1+\delta)\kappa^{1+\delta}}
		\int_Jw(x)^{1+\delta}\,dx.
		\]
		All these integrals are finite because $w\in L^\infty(U)$
		and $U$ is bounded.
		
		Choose $\delta>0$ sufficiently small that
		\[
		\frac{2\delta}
		{(1+\delta)\kappa^{1+\delta}}
		\leq\frac12.
		\]
		Setting $r=1+\delta$, we conclude that
		\begin{equation}\label{weight-weak-RH-local}
			\int_J w(x)^r\,dx
			\leq
			2\bigl[M^-(w\chi_J)(b)\bigr]^{r-1}w(J)
		\end{equation}
		for every $J=(a,b)\subset U$.
		The exponent $r$ depends only on $\beta_{1/4}$.
		
		It remains to obtain \eqref{local-weight-RH}.
		Fix $h>0$ and $a$ such that
		\[
		u\leq a-h<a<a+h\leq v.
		\]
		Write
		\[
		\ell=a-h,\qquad m=a,\qquad z=a+h,
		\]
		and set
		\[
		A=\int_\ell^m w(x)^r\,dx,
		\qquad
		W=w(\ell,z).
		\]
		For every $x\in(m,z)$, apply
		\eqref{weight-weak-RH-local} to $(\ell,x)$ to obtain
		\begin{align*}
			A
			&\leq\int_\ell^x w(y)^r\,dy\\
			&\leq
			2\bigl[M^-(w\chi_{(\ell,x)})(x)\bigr]^{r-1}w(\ell,x)\\
			&\leq
			2\bigl[M^-(w\chi_{(\ell,z)})(x)\bigr]^{r-1}W.
		\end{align*}
		Hence
		\[
		M^-(w\chi_{(\ell,z)})(x)
		\geq
		\left(\frac{A}{2W}\right)^{1/(r-1)}
		\qquad(m<x<z).
		\]
		The weak type $(1,1)$ estimate for $M^-$, applied
		at half this threshold, yields
		\[
		h
		\leq
		C W\left(\frac{2W}{A}\right)^{1/(r-1)}.
		\]
		Rearranging gives
		\begin{equation}\label{weight-total-mass-local}
			A\leq C_r h^{1-r}W^r.
		\end{equation}
		
		Finally, apply \eqref{local-weight-condition}
		with the triple $\ell<m<z$, the set $E=(m,z)$,
		and $\eta=1/2$.
		Since
		\[
		\frac{|E|}{z-\ell}=\frac12,
		\]
		the conclusion of that implication is false.
		Consequently $w(m,z)>\beta_{1/2}w(\ell,m)$ and therefore,
		\[
		W\leq
		(1+\beta_{1/2}^{-1})w(m,z).
		\]
		Combining this estimate with
		\eqref{weight-total-mass-local}, we obtain
		\[
		\left(\frac1h\int_\ell^m w(x)^r\,dx\right)^{1/r}
		\leq
		C\frac1h\int_m^z w(x)\,dx,
		\]
		which is \eqref{local-weight-RH}.
		The constants depend only on $\beta_{1/4}$
		and $\beta_{1/2}$, as asserted.
	\end{proof}
	
	The next lemma is the one-sided analogue of
	\cite[Lemma~4.56]{CruzUribeFiorenza}.
	In the two-sided setting, the proof first derives a uniform
	$A_\infty$-type property for the family of weights
	\[
	w_t(x)=t^{p(x)}
	\]
	from the modular estimate of \cite[Lemma~4.55]{CruzUribeFiorenza}, and then applies a reverse
	H\"older self-improvement argument.
	
	In the present one-sided setting, we combine
	Lemma~\ref{lemma455-CU} with the classical theory of the classes
	$A_\infty^+$ and $RH_s^+$ developed in
	\cite{MartinReyesPickTorre1993,CruzUribe1997,RiverosTorre2001}.
	
	\begin{lemma}\label{lemma456-CU}
		Let $p(\cdot)\in K_0^+(\mathbb R)$ satisfy $1<p_-\leq p_+<\infty$, and fix the constants $b_1,b_2>0$.
		Let
		\[
		I_0^-=(a_0-h_0,a_0),
		\qquad
		I_0^+=(a_0,a_0+h_0),
		\]
		where $h_0>0$, and set $U=(a_0-h_0,a_0+h_0)$.
		Suppose that $b_1\leq t\leq \frac{b_2}{\|\chi_U\|_{p(\cdot)}}.$
		Then there exist $s_0>1$ and $C_0>0$, depending only on
		$p_-,p_+,b_1,b_2$, and $[p]_{K_0^+}$, such that
		\begin{equation}\label{456plus}
			\left(
			\fint_{I_0^-}t^{sp(x)}\,dx
			\right)^{1/s}
			\leq
			C_0\fint_{I_0^+}t^{p(x)}\,dx
		\end{equation}
		for every $1<s<s_0$.
	\end{lemma}
	
	\begin{proof}
		Set $w(x)=t^{p(x)}$. We first verify
		condition~\eqref{local-weight-condition} of
		Lemma~\ref{local-one-sided-RH}, with constants
		independent of $t$ and $U$. Fix $\eta\in(0,1)$.
		We claim that there exists $\beta_\eta\in(0,1)$,
		independent of $t$ and $U$, such that
		\begin{equation}\label{local-Ainf-U}
			w(E)\leq\beta_\eta w(a,b)
			\quad\Longrightarrow\quad
			\frac{|E|}{c-a}<\eta
		\end{equation}
		whenever $a_0-h_0\leq a<b<c\leq a_0+h_0$ and $E\subset(b,c)$ is measurable.
		
		Indeed, since $(a,b)\subset U$, our hypothesis for $t$ gives
		\[
		\|t\chi_{(a,b)}\|_{p(\cdot)}
		\leq b_2.
		\]
		Consequently $w(a,b)\leq B, \ B:=\max\{b_2^{p_-},b_2^{p_+}\}.$
		
		Suppose that $0<\beta\leq1$ and $w(E)\leq\beta w(a,b).$
		Then $\rho_{p(\cdot)}(t\chi_E)=w(E)\leq B,$
		and hence
		\[
		\|t\chi_E\|_{p(\cdot)}
		\leq c_2,
		\qquad
		c_2:=\max\{1,B^{1/p_-}\}.
		\]
		Suppose, on the contrary, that $|E|\geq\eta(c-a),$
		then
		\[
		\frac{1}{c-a}\int_b^c t\chi_E(y)\,dy
		\geq b_1\eta.
		\]
		Lemma~\ref{lemma455-CU}, applied to $t\chi_E$, yields
		\[
		\int_a^b
		\left(t\frac{|E|}{c-a}\right)^{p(x)}\,dx
		\leq C_\eta w(E),
		\]
		where $C_\eta$ is independent of $t,U,a,b,c$, and $E$.
		On the other hand,
		\[
		\int_a^b
		\left(t\frac{|E|}{c-a}\right)^{p(x)}\,dx
		\geq \eta^{p_+}w(a,b).
		\]
		Choosing
		\[
		\beta_\eta
		=
		\min\left\{
		\frac12,\frac{\eta^{p_+}}{2C_\eta}
		\right\}
		\]
		gives a contradiction.
		This proves \eqref{local-Ainf-U}.
		
		Since $t>0$ is fixed and $p_+<\infty$, the weight
		$w(x)=t^{p(x)}$ is positive almost everywhere and belongs
		to $L^\infty(U)$.
		Moreover, \eqref{local-Ainf-U} is precisely the
		hypothesis of Lemma~\ref{local-one-sided-RH}, with
		constants independent of $t$ and $U$.
		Therefore, there exist $r>1$ and $C>0$, independent
		of $t$ and $U$, such that
		\[
		\left(
		\fint_{I_0^-}t^{rp(x)}\,dx
		\right)^{1/r}
		\leq
		C\fint_{I_0^+}t^{p(x)}\,dx.
		\]
		For every $1<s<r$, H\"older's inequality gives
		\[
		\left(
		\fint_{I_0^-}t^{sp(x)}\,dx
		\right)^{1/s}
		\leq
		\left(
		\fint_{I_0^-}t^{rp(x)}\,dx
		\right)^{1/r}.
		\]
		Taking $s_0=r$ proves \eqref{456plus}, with a constant
		independent of $s$.
	\end{proof}
	
	The following lemma is a consequence of the boundedness
	criterion of Lerner and Ombrosi
	\cite[Theorem~1.2]{LernerOmbrosi2010}.
	We include the proof to explain its application to the right-sided maximal operator.
	
	\begin{lemma}\label{passage-to-p}
		Let $p(\cdot)\in\mathcal P(\mathbb R)$ satisfy
		$1<p_-\leq p_+<\infty$.
		Suppose that there exist constants $s_*>1$ and $C>0$
		such that
		\begin{equation}\label{passage-hypothesis}
			\|M^+\|_{L^{sp(\cdot)}(\mathbb R)
				\to L^{sp(\cdot)}(\mathbb R)}
			\leq C(s-1)^{-1/(sp_-)}
			\qquad (1<s<s_*),
		\end{equation}
		where $C$ is independent of $s$.
		Then $M^+$ is bounded on $L^{p(\cdot)}(\mathbb R)$.
	\end{lemma}
	
	\begin{proof}
		For $0<r<1$, define
		\[
		M_r^+f
		=
		\bigl(M^+(|f|^r)\bigr)^{1/r}.
		\]
		The basis $\{(x,x+h)\}_{h>0}$, $x\in\mathbb R$, which generates $M^+$,
		satisfies Stein's property \cite[Definition~1.1]{LernerOmbrosi2010};
		this was already observed by Sawyer \cite{Sawyer1986}, as noted in
		\cite[p.~55]{LernerOmbrosi2010}. Moreover, $L^{p(\cdot)}(\mathbb R)$ is a Banach function space. Therefore, by
		\cite[Theorem~1.2, $(i)\Rightarrow(ii)$]
		{LernerOmbrosi2010}, it suffices to prove that
		\begin{equation}\label{LO-condition}
			\lim_{\varepsilon\to0^+}
			\varepsilon
			\|M_{1-\varepsilon}^+\|_
			{L^{p(\cdot)}\to L^{p(\cdot)}}=0.
		\end{equation}
		
		Set $r=1-\varepsilon$ and $s=1/r$.
		By the Lemma \ref{dilationExpo} we have,
		\[
		\|M_r^+f\|_{p(\cdot)}
		=
		\|M^+(|f|^r)\|_{sp(\cdot)}^s,
		\qquad
		\||f|^r\|_{sp(\cdot)}
		=
		\|f\|_{p(\cdot)}^r.
		\]
		Since every non-negative function in $L^{sp(\cdot)}$
		can be written as $|f|^r$ with
		$f\in L^{p(\cdot)}$, it follows that
		\begin{equation}\label{Mr-Msp}
			\|M_r^+\|_{L^{p(\cdot)}\to L^{p(\cdot)}}
			=
			\left(
			\|M^+\|_{L^{sp(\cdot)}\to L^{sp(\cdot)}}
			\right)^s.
		\end{equation}
		Using $\varepsilon=(s-1)/s$ and
		\eqref{passage-hypothesis}, we obtain
		\[
		\varepsilon
		\|M_{1-\varepsilon}^+\|_
		{L^{p(\cdot)}\to L^{p(\cdot)}}
		\leq
		\frac{C^s}{s}(s-1)^{1-1/p_-}.
		\]
		The right-hand side tends to zero as $s\to1^+$,
		because $p_->1$.
		Thus \eqref{LO-condition} holds, and the
		Lerner--Ombrosi criterion yields the conclusion.
	\end{proof}
	
	\begin{proof}[Proof of Theorem~\ref{main-Theorem}]
		We are now ready to prove Theorem~\ref{main-Theorem}.
		The argument follows the general strategy of
		\cite[Theorem~4.52]{CruzUribeFiorenza}.
		In place of the usual Calder\'on--Zygmund decomposition,
		we use the one-sided level-set decomposition for $M^+$
		employed in
		\cite[p.~525, proof of Theorem~2, equation~(5.1)]
		{MartinReyesOrtegaTorre1990}, with $g\equiv1$.
		We then combine Lemmas~\ref{lemma455-CU}
		and~\ref{lemma456-CU} on consecutive intervals
		adapted to this decomposition.
		
		We shall work with $1<s<s_0$, where $s_0>1$ is furnished
		by Lemma~\ref{lemma456-CU} with $b_1=1/2$ and the fixed
		constant $b_2$ determined in \eqref{tn-upper}.
		These constants are independent of $s$. Replacing $s_0$ by $\min\{s_0,2\}$, we may assume that $s_0\leq2$. As in the proof of
		\cite[Theorem~4.52]{CruzUribeFiorenza}, we may reduce to the case in
		which $f$ is non-negative, bounded, compactly supported, and
		\[
		\|f\|_{sp(\cdot)}=1.
		\]
		We decompose $f=f_1+f_2,$
		\[
		f_1=f\chi_{\{f>1\}},
		\qquad
		f_2=f\chi_{\{f\leq1\}}.
		\]
		Since $f_1>1$ on its support and $s>1$,
		\[
		\|f_1\|_{p(\cdot)}
		\leq
		\|f_1\|_{sp(\cdot)}
		\leq1.
		\]
		
		We first estimate the large part $f_1$. Fix $A>1$ and, for every
		integer $k\geq0$, define
		\[
		\Omega_k^+
		=
		\left\{
		x\in\mathbb R:
		M^+f_1(x)>A^k
		\right\}.
		\]
		By this one-sided level-set decomposition,
		\[
		\Omega_k^+
		=
		\bigcup_j I_j^k,
		\qquad
		I_j^k=(a_j^k,b_j^k),
		\]
		where the intervals $I_j^k$ are pairwise disjoint and
		\begin{equation}\label{rising-average}
			\frac{1}{|I_j^k|}
			\int_{I_j^k}f_1(y)\,dy
			=
			A^k.
		\end{equation}
		Moreover, if $x\in I_j^k$, then
		\begin{equation}\label{rising-tail}
			\frac{1}{b_j^k-x}
			\int_x^{b_j^k}f_1(y)\,dy
			>
			A^k.
		\end{equation}
		This is the one-dimensional one-sided substitute for the
		Calder\'on--Zygmund cubes used in the proof of the two-sided result;
		see, for instance,
		\cite{Sawyer1986,MartinReyesOrtegaTorre1990,MartinReyes1993}.
		
		We shall also use the level pieces
		\[
		E_j^k
		=
		I_j^k\cap
		\left\{
		M^+f_1\leq A^{k+1}
		\right\}.
		\]
		Since
		\[
		\Omega_k^+\setminus\Omega_{k+1}^+
		=
		\bigcup_jE_j^k,
		\]
		the sets $E_j^k$ are pairwise disjoint as $j$ and $k$ vary. In addition, the level-set decomposition gives
		\begin{equation}\label{level-sparse}
			|\Omega_{k+1}^+\cap I_j^k|
			\leq
			\frac1A|I_j^k|.
		\end{equation}
		Indeed,
		\begin{align*}
			A^{k+1}|\Omega_{k+1}^+\cap I_j^k|
			&=
			\sum_{I_i^{k+1}\subset I_j^k}
			A^{k+1}|I_i^{k+1}|
			\\
			&=
			\sum_{I_i^{k+1}\subset I_j^k}
			\int_{I_i^{k+1}}f_1(y)\,dy
			\\
			&\leq
			\int_{I_j^k}f_1(y)\,dy
			=
			A^k|I_j^k|.
		\end{align*}
		Consequently,
		\begin{equation}\label{Ej-large}
			|E_j^k|
			\geq
			\left(1-\frac1A\right)|I_j^k|.
		\end{equation}
		This is the one-sided analogue of the large-subset property used in
		the classical proof.
		
		We now work inside a fixed component
		\[
		I=I_j^k=(a,b),
		\qquad
		\lambda=A^k.
		\]
		The following geometric partition is chosen so that
		Lemmas~\ref{lemma456-CU} and \ref{lemma455-CU} can be applied
		successively. Let
		\[
		\ell_0=\frac{b-a}{3},
		\qquad
		\ell_{n+1}=\frac23\ell_n,
		\qquad n\geq0,
		\]
		and define
		\[
		L_n=(b-3\ell_n,b-2\ell_n),
		\qquad
		M_n=(b-2\ell_n,b-\ell_n),
		\qquad
		R_n=(b-\ell_n,b).
		\]
		For each $n$, these are three consecutive intervals of equal length
		$\ell_n$. The intervals $\{L_n\}_{n\geq0}$ are pairwise disjoint and
		\[
		I=\bigcup_{n=0}^\infty L_n
		\]
		up to a set of measure zero.
		
		Set
		\begin{equation}\label{tn-definition}
			t_n
			=
			\frac{1}{2\ell_n}
			\int_{R_n}f_1(y)\,dy.
		\end{equation}
		By \eqref{rising-tail}, applied at the left endpoint of $R_n$,
		\[
		\frac1{\ell_n}\int_{R_n}f_1(y)\,dy>\lambda,
		\]
		and therefore
		\begin{equation}\label{tn-lower}
			t_n>\frac{\lambda}{2}\geq\frac12.
		\end{equation}
		
		We next verify the upper bound required in
		Lemma~\ref{lemma456-CU}.
		By H\"older's inequality and
		$\|f_1\|_{p(\cdot)}\leq1$,
		\[
		t_n
		=
		\frac{1}{2\ell_n}\int_{R_n}f_1(y)\,dy
		\leq
		\frac{C}{2\ell_n}
		\|\chi_{R_n}\|_{p'(\cdot)}.
		\]
		Apply the $K_0^+$ condition to the triple
		\[
		b-3\ell_n<b-\ell_n<b.
		\]
		The corresponding left interval agrees, up to a set of
		measure zero, with $L_n\cup M_n$, while the right interval
		is $R_n$. Hence
		\[
		\|\chi_{L_n\cup M_n}\|_{p(\cdot)}
		\|\chi_{R_n}\|_{p'(\cdot)}
		\leq
		3[p]_{K_0^+}\ell_n.
		\]
		Consequently,
		\begin{equation}\label{tn-upper}
			t_n
			\leq
			\frac{C}{\|\chi_{L_n\cup M_n}\|_{p(\cdot)}}.
		\end{equation}
		Together with \eqref{tn-lower}, this verifies the
		hypotheses of Lemma~\ref{lemma456-CU} for
		\[
		I_0^-=L_n,\qquad I_0^+=M_n,
		\]
		with $b_1=1/2$ and a fixed constant $b_2>0$ independent
		of $j,k,n$, and $s$.
		
		Applying Lemma~\ref{lemma456-CU} to the adjacent intervals $L_n$
		and $M_n$, we obtain
		\begin{equation}\label{RH-local}
			\left(
			\fint_{L_n}t_n^{sp(x)}\,dx
			\right)^{1/s}
			\leq
			C
			\fint_{M_n}t_n^{p(x)}\,dx.
		\end{equation}
		On the other hand, from the definition of $t_n$,
		\[
		t_n
		=
		\frac{1}{|M_n|+|R_n|}
		\int_{R_n}f_1(y)\,dy.
		\]
		Therefore Lemma~\ref{lemma455-CU}, applied to the triple
		\[
		b-2\ell_n<b-\ell_n<b,
		\]
		gives
		\begin{equation}\label{modular-local}
			\int_{M_n}t_n^{p(x)}\,dx
			\leq
			C
			\int_{R_n}f_1(y)^{p(y)}\,dy.
		\end{equation}
		Since
		\[
		|L_n|=|M_n|=|R_n|=\ell_n,
		\]
		\eqref{RH-local} and \eqref{modular-local} imply
		\begin{equation}\label{local-combined}
			\int_{L_n}t_n^{sp(x)}\,dx
			\leq
			C\ell_n
			\left(
			\fint_{R_n}f_1(y)^{p(y)}\,dy
			\right)^s.
		\end{equation}
		Using \eqref{tn-lower}, we conclude that
		\begin{equation}\label{lambda-local}
			\int_{L_n}\lambda^{sp(x)}\,dx
			\leq
			C\ell_n
			\left(
			\fint_{R_n}f_1(y)^{p(y)}\,dy
			\right)^s.
		\end{equation}
		Thus \eqref{lambda-local} is the one-sided counterpart of the
		combined use of the modular and reverse H\"older estimates in the
		proof of \cite[Theorem~4.52]{CruzUribeFiorenza}.
		
		It remains to sum these local inequalities over the level intervals.
		The intervals $I_j^k$ are disjoint for fixed $k$, but they may be
		nested as $k$ varies. To handle this point, we use Lemma~\ref{sparse-alike}.
		
		For each $k\geq0$, let
		\[
		F_k
		=
		\left\{
		x\in\mathbb R:
		M^+f_1(x)\leq A^{k+2}
		\right\}.
		\]
		Fix $n$ such that $E_j^k\cap L_n\neq\varnothing$, and choose
		$x_n^*\in E_j^k\cap L_n$. Since $M^+f_1(x_n^*)\leq A^{k+1},$ we can apply Lemma~\ref{sparse-alike} with $v=f_1$, $x=x_n^*$, $y=b-\ell_n$, \ $\lambda_1=A^{k+1}, \ \lambda_2=A^{k+2}.$
		
		We define
		\[
		J_n=(x_n^*,b-\ell_n).
		\]
		Because $x_n^*\in L_n$, $|J_n|\geq\ell_n$ and, hence,
		\begin{equation}\label{Jn-density}
			|F_k\cap J_n|
			\geq
			\left(1-\frac1A\right)|J_n|
			\geq
			\left(1-\frac1A\right)\ell_n.
		\end{equation}
		
		Set $g(y)=f_1(y)^{p(y)}.$ If $z\in J_n$, then $(z,b)$ contains $R_n$ and $b-z\leq3\ell_n.$
		Therefore
		\[
		M^+g(z)
		\geq
		\frac{1}{b-z}\int_{R_n}g(y)\,dy
		\geq
		\frac13\fint_{R_n}g(y)\,dy.
		\]
		Using \eqref{lambda-local} and \eqref{Jn-density}, we obtain
		\begin{align}
			\int_{E_j^k\cap L_n}\lambda^{sp(x)}\,dx
			&\leq
			C\ell_n
			\left(
			\fint_{R_n}g(y)\,dy
			\right)^s
			\nonumber\\
			&\leq
			C
			\int_{F_k\cap J_n}
			\bigl(M^+g(z)\bigr)^s\,dz.
			\label{local-sparse-sum}
		\end{align}
		
		The intervals $J_n$ have uniformly bounded overlap. Indeed,
		\[
		J_n\subset(b-3\ell_n,b-\ell_n),
		\]
		and since $\ell_{n+1}=\frac23\ell_n$, any point belongs to at most a
		fixed number of these intervals. Thus
		\[
		\sum_n\chi_{J_n}\leq C.
		\]
		Summing \eqref{local-sparse-sum} in $n$ gives
		\begin{equation}\label{sum-fixed-jk}
			\int_{E_j^k}(A^k)^{sp(x)}\,dx
			\leq
			C
			\int_{F_k\cap I_j^k}
			\bigl(M^+g(z)\bigr)^s\,dz.
		\end{equation}
		
		For fixed $k$, the intervals $I_j^k$ are pairwise disjoint.
		Consequently,
		\[
		\sum_j
		\int_{E_j^k}(A^k)^{sp(x)}\,dx
		\leq
		C
		\int_{\Omega_k^+\cap F_k}
		\bigl(M^+g(z)\bigr)^s\,dz.
		\]
		But
		\[
		\Omega_k^+\cap F_k
		\subset
		\left\{
		A^k<M^+f_1\leq A^{k+2}
		\right\}.
		\]
		These bands have uniformly bounded overlap as $k$ varies. Therefore,
		\begin{equation}\label{global-large-sum}
			\sum_{k,j}
			\int_{E_j^k}(A^k)^{sp(x)}\,dx
			\leq
			C
			\int_{\mathbb R}
			\left[
			M^+\!\left(f_1^{p(\cdot)}\right)(x)
			\right]^s
			dx.
		\end{equation}
		
		We now use the classical $L^s$ estimate for the one-sided maximal
		operator with its precise dependence on $s$. For $1<s\leq2$,
		\begin{equation}\label{classical-Mplus-Ls}
			\int_{\mathbb R}(M^+g(x))^s\,dx
			\leq
			\frac{C}{s-1}
			\int_{\mathbb R}|g(x)|^s\,dx,
		\end{equation}
		where $C$ is independent of $s$.
		
		For completeness, this estimate follows directly from the weak type
		$(1,1)$ inequality. Indeed, for every $\lambda>0$,
		\[
		M^+\!\left(g\chi_{\{|g|\leq\lambda/2\}}\right)
		\leq\frac{\lambda}{2},
		\]
		and hence
		\[
		\{M^+g>\lambda\}
		\subset
		\left\{
		M^+\!\left(g\chi_{\{|g|>\lambda/2\}}\right)
		>\frac{\lambda}{2}
		\right\}.
		\]
		The weak type $(1,1)$ estimate therefore gives
		\[
		|\{M^+g>\lambda\}|
		\leq
		\frac{C}{\lambda}
		\int_{\{|g|>\lambda/2\}}|g(y)|\,dy.
		\]
		Using the distribution-function formula and Fubini's theorem,
		\begin{align*}
			\int_{\mathbb R}(M^+g(x))^s\,dx
			&=
			s\int_0^\infty
			\lambda^{s-1}
			|\{M^+g>\lambda\}|\,d\lambda
			\\
			&\leq
			Cs\int_{\mathbb R}|g(y)|
			\int_0^{2|g(y)|}\lambda^{s-2}\,d\lambda\,dy
			\\
			&\leq
			\frac{C}{s-1}
			\int_{\mathbb R}|g(y)|^s\,dy,
		\end{align*}
		where in the last inequality we used $1<s\leq2$.
		
		Applying this estimate to
		\[
		g=f_1^{p(\cdot)},
		\]
		we obtain
		\[
		\int_{\mathbb R}
		\left[
		M^+\!\left(f_1^{p(\cdot)}\right)(x)
		\right]^s dx
		\leq
		\frac{C}{s-1}
		\int_{\mathbb R}f_1(x)^{sp(x)}\,dx
		\leq
		\frac{C}{s-1}.
		\]
		
		Since the sets $E_j^k$ decompose $\{M^+f_1>1\}$ and $M^+f_1(x)\leq A^{k+1}$ on $E_j^k,$ we have
		\[
		(M^+f_1(x))^{sp(x)}
		\leq
		A^{sp_+}(A^k)^{sp(x)}
		\qquad\text{on }E_j^k.
		\]
		Since $s_0\leq2$, the factor $A^{sp_+}$ is uniformly bounded. Therefore,
		\begin{equation}\label{large-part-final}
			\int_{\{M^+f_1>1\}}
			\left(M^+f_1(x)\right)^{sp(x)}\,dx
			\leq
			\frac{C}{s-1}.
		\end{equation}
		
		The complementary region is simpler. On $\{M^+f_1\leq1\}$ we have,
		since $p(x)\geq p_-$,
		\[
		\left(M^+f_1(x)\right)^{sp(x)}
		\leq
		\left(M^+f_1(x)\right)^{sp_-}.
		\]
		Since $sp_->1$, the classical $L^{sp_-}$ boundedness of $M^+$ and
		the fact that $f_1>1$ on its support yield
		\[
		\int_{\{M^+f_1\leq1\}}
		\left(M^+f_1(x)\right)^{sp(x)}\,dx
		\leq
		C
		\int_{\mathbb R}f_1(x)^{sp_-}\,dx
		\leq
		C
		\int_{\mathbb R}f_1(x)^{sp(x)}\,dx
		\leq C.
		\]
		Combining this estimate with \eqref{large-part-final}, we obtain
		\begin{equation}\label{large-part-modular}
			\int_{\mathbb R}
			\left(M^+f_1(x)\right)^{sp(x)}\,dx
			\leq
			\frac{C}{s-1}.
		\end{equation}
		
		We next turn to the small part $f_2$. Since $M^+f_2\le Mf_2$, it suffices
		to bound $\|Mf_2\|_{sp(\cdot)}$ by a constant independent of $s$. We cannot
		invoke the boundedness of $M$ on $L^{sp(\cdot)}(\mathbb R)$, since
		$p(\cdot)$ need not belong to $K_0$ (see Example~\ref{strict-one-sided}).
		However, the estimate of $Mf_2$ in the proof of
		\cite[Theorem~4.52]{CruzUribeFiorenza} uses only the $N_\infty$ condition:
		it is the argument given in the proof of
		\cite[Theorem~3.16, pp.~104--106]{CruzUribeFiorenza}, where, as observed in
		\cite[pp.~131--132]{CruzUribeFiorenza}, the decay of the exponent at
		infinity enters only through the embeddings of $L^\infty$ into the defect
		spaces associated with $p(\cdot)$ and $p_\infty$, and these embeddings
		follow from $N_\infty$. This is essentially the mechanism introduced by
		Nekvinda~\cite{Nekvinda2004}.
		
		Since $s>1$,
		\[
		\left|\frac{1}{sp(x)}-\frac{1}{sp_\infty}\right|^{-1}
		= s\left|\frac{1}{p(x)}-\frac{1}{p_\infty}\right|^{-1}
		\ge \left|\frac{1}{p(x)}-\frac{1}{p_\infty}\right|^{-1},
		\]
		so $sp(\cdot)\in N_\infty$ with limiting exponent $sp_\infty$ and the same
		constant $\Lambda_\infty$. Moreover, $p_\infty\in[p_-,p_+]$ and $1<s<s_0\le2$,
		so $sp_\infty$ ranges over the compact interval $[p_-,2p_+]\subset(1,\infty)$,
		and the norm of $M$ on $L^{sp_\infty}(\mathbb R)$ is bounded uniformly in $s$.
		Therefore, exactly as in \cite[p.~169]{CruzUribeFiorenza},
		\begin{equation}\label{small-part-final}
			\|M^+f_2\|_{sp(\cdot)}
			\le
			\|Mf_2\|_{sp(\cdot)}
			\le C,
		\end{equation}
		where $C$ depends only on $p(\cdot)$ and not on $s$.
		
		
		We now combine the estimates for $f_1$ and $f_2$.
		By \eqref{large-part-modular} and the
		modular--norm inequalities,
		\[
		\|M^+f_1\|_{sp(\cdot)}
		\leq C(s-1)^{-1/(sp_-)}.
		\]
		Indeed, after enlarging $C$ and restricting to
		$1<s<s_*\leq2$, the modular bound is at least one,
		and the lower exponent of $sp(\cdot)$ is $sp_-$.
		On the other hand, \eqref{small-part-final} gives
		\[
		\|M^+f_2\|_{sp(\cdot)}\leq C.
		\]
		Here $s_*>1$ is chosen sufficiently close to one
		so that both estimates hold with constants
		independent of $s$.
		
		By sublinearity, we therefore have
		\[
		\|M^+f\|_{sp(\cdot)}
		\leq C(s-1)^{-1/(sp_-)}
		\]
		for every non-negative, bounded, compactly
		supported function $f$ with
		$\|f\|_{sp(\cdot)}=1$.
		Homogeneity extends this estimate to all bounded,
		compactly supported functions.
		
		For an arbitrary $f\in L^{sp(\cdot)}(\mathbb R)$,
		apply the estimate to
		\[
		f_N=\min\{|f|,N\}\chi_{(-N,N)}.
		\]
		Then $f_N\uparrow|f|$ pointwise, and hence, by the monotone convergence
		theorem applied to each average, $M^+f_N\uparrow M^+f$ pointwise.
		Since $\|f_N\|_{sp(\cdot)}\le\|f\|_{sp(\cdot)}$, Fatou's lemma in variable
		Lebesgue spaces \cite[Theorem~2.61]{CruzUribeFiorenza} yields
		\begin{equation}\label{Mplus-sp-norm}
			\|M^+\|_{L^{sp(\cdot)}\to L^{sp(\cdot)}}
			\leq C(s-1)^{-1/(sp_-)}
			\qquad (1<s<s_*).
		\end{equation}
		Lemma~\ref{passage-to-p} now implies that
		$M^+$ is bounded on $L^{p(\cdot)}(\mathbb R)$.
		This completes the proof of
		Theorem~\ref{main-Theorem}.
	\end{proof}
	
	\section{Examples}
	\label{examples}
	
	We present three examples illustrating the relation between
	$K_0$, $\Kplus$, and $\Ninf$.
	Recall that the classical condition $K_0$ is defined by
	\[
	[p]_{K_0}
	:=
	\sup_I
	\frac{
		\|\chi_I\|_{p(\cdot)}
		\|\chi_I\|_{p'(\cdot)}
	}{|I|}
	<\infty,
	\]
	where the supremum is taken over all bounded intervals.
	By monotonicity of the Luxemburg norm, whenever $a<b<c$,
	\[
	\|\chi_{(a,b)}\|_{p(\cdot)}
	\|\chi_{(b,c)}\|_{p'(\cdot)}
	\le
	\|\chi_{(a,c)}\|_{p(\cdot)}
	\|\chi_{(a,c)}\|_{p'(\cdot)}.
	\]
	Consequently,
	\begin{equation}\label{examples-K0-inclusion}
		K_0\subset\Kplus,
		\qquad
		[p]_{\Kplus}\le [p]_{K_0}.
	\end{equation}
	
	We shall repeatedly use the elementary estimate
	\begin{equation}\label{examples-characteristic}
		\|\chi_E\|_{q(\cdot)}
		\le |E|^{1/q_+(E)},
		\qquad 0<|E|\le1,
	\end{equation}
	where $q_+(E)=\operatorname*{ess\,sup}_{E}q$.
	Indeed, the modular of
	$\chi_E/|E|^{1/q_+(E)}$ is at most one.
	In particular, if $I$ is an interval with $|I|\le1$, then
	\begin{equation}\label{examples-small-interval}
		\frac{
			\|\chi_I\|_{p(\cdot)}
			\|\chi_I\|_{p'(\cdot)}
		}{|I|}
		\le
		|I|^{1/p_+(I)-1/p_-(I)}.
	\end{equation}
	
	\begin{example}[The condition $\Kplus$ does not imply $\Ninf$]
		\label{Kplus-without-Ninf}
		Let $\varphi\in C_c^\infty(\R)$ satisfy
		\[
		0\le\varphi\le1,\qquad
		\operatorname{supp}\varphi\subset[-1/2,1/2],
		\qquad
		\varphi=1\quad\text{on }[-1/4,1/4].
		\]
		Define
		\[
		p(x)=3+\sum_{k=1}^{\infty}\varphi(x-e^k).
		\]
		Then
		\[
		p(\cdot) \in K_0\subset\Kplus,
		\qquad
		p(\cdot)\notin\Ninf.
		\]
		This is the exponent from
		\cite[Example~4.51]{CruzUribeFiorenza}.
	\end{example}
	
	\begin{proof}
		The supports of the summands are pairwise disjoint.
		Thus $3\le p\le4$, and $p$ is globally Lipschitz.
		Let $L$ be a Lipschitz constant.
		
		Fix an interval $I$ and write $D=|I|$.
		If $D\le1$, then
		\[
		p_+(I)-p_-(I)\le LD,
		\qquad
		\frac1{p_-(I)}-\frac1{p_+(I)}
		\le\frac{LD}{9}.
		\]
		Therefore, by \eqref{examples-small-interval},
		\[
		\frac{
			\|\chi_I\|_{p(\cdot)}
			\|\chi_I\|_{p'(\cdot)}
		}{D}
		\le D^{-LD/9}
		\le \exp\!\left(\frac{L}{9e}\right).
		\]
		
		Suppose now that $D\ge1$, and set
		\[
		S=\bigcup_{k\ge1}[e^k-1/2,e^k+1/2].
		\]
		There is an absolute constant $C$ such that
		\begin{equation}\label{examples-sparse-support}
			|I\cap S|\le C\log(e+D).
		\end{equation}
		Indeed, if the first and last intervals in this union
		meeting $I$ have indices $m$ and $n$, then
		$e^n-e^m\le D+1$, whence
		\[
		n-m+1
		\le 1+\log\!\left(1+\frac{D+1}{e}\right).
		\]
		Each of these intervals has length one.
		
		Since $p\ge3$ and $D\ge1$,
		\[
		\int_I \bigl(D^{-1/3}\bigr)^{p(x)}\,dx\le1,
		\qquad
		\|\chi_I\|_{p(\cdot)}\le D^{1/3}.
		\]
		Moreover, $4/3\le p'\le3/2$ and $p'=3/2$ outside $S$.
		For $A\ge1$, \eqref{examples-sparse-support} gives
		\[
		\int_I
		\left(\frac1{AD^{2/3}}\right)^{p'(x)}\,dx
		\le
		A^{-3/2}
		+C A^{-4/3}D^{-8/9}\log(e+D).
		\]
		The function $D^{-8/9}\log(e+D)$ is bounded for $D\ge1$.
		Choosing $A$ sufficiently large, independently of $I$,
		we obtain
		\[
		\|\chi_I\|_{p'(\cdot)}\le AD^{2/3}.
		\]
		This proves $p\in K_0$, and hence $p\in\Kplus$.
		
		Finally, suppose that $p\in\Ninf$ with limiting exponent
		$p_\infty$.
		Since $p=3$ on $(-\infty,0)$, the defining integrability
		condition forces $p_\infty=3$.
		However, $p=4$ on
		\[
		E=\bigcup_{k\ge1}[e^k-1/4,e^k+1/4],
		\qquad |E|=\infty.
		\]
		On this set,
		\[
		\left|\frac1{p(x)}-\frac13\right|=\frac1{12}.
		\]
		Thus, for every $\Lambda>0$,
		\[
		\int_E
		\exp\!\left(
		-\Lambda
		\left|\frac1{p(x)}-\frac13\right|^{-1}
		\right)\,dx
		=
		\int_E e^{-12\Lambda}\,dx
		=\infty.
		\]
		Consequently, $p\notin\Ninf$.
	\end{proof}
	
	\begin{example}[An exponent without local log-H\"older continuity]
		\label{non-log-holder}
		Fix $0<\alpha<1$ and let
		\[
		r_\alpha=\exp(-3^{1/\alpha}).
		\]
		Define
		\[
		p_\alpha(x)=
		\begin{cases}
			\displaystyle
			\left(\frac12+
			\bigl(\log(1/|x|)\bigr)^{-\alpha}\right)^{-1},
			&0<|x|<r_\alpha,\\[6pt]
			2,    &x=0,\\[2pt]
			\displaystyle\frac65,
			&|x|\ge r_\alpha.
		\end{cases}
		\]
		Then
		\[
		p_\alpha\in K_0\cap\Ninf
		\subset\Kplus\cap\Ninf,
		\]
		but $p_\alpha$ is not locally log-H\"older continuous
		at the origin.
		This construction is taken from
		\cite[Example~4.59]{CruzUribeFiorenza},
		with the constant extension matching the boundary value.
	\end{example}
	
	\begin{proof}
		Write $p=p_\alpha$ and $r=r_\alpha$.
		The definition gives
		\[
		\frac65\le p(x)\le2.
		\]
		Since $p=6/5$ outside $(-r,r)$, the condition $\Ninf$
		holds with $p_\infty=6/5$: its defining integrand is
		bounded by one and is supported on a set of finite measure.
		
		On the other hand, for $0<|x|<r$,
		\[
		2-p(x)
		=
		\frac{
			4\bigl(\log(1/|x|)\bigr)^{-\alpha}
		}{
			1+2\bigl(\log(1/|x|)\bigr)^{-\alpha}
		}.
		\]
		Hence
		\[
		|p(x)-p(0)|\log(1/|x|)
		\longrightarrow\infty
		\qquad\text{as }x\to0.
		\]
		Thus local log-H\"older continuity fails.
		
		It remains to verify $K_0$.
		We first consider intervals contained in $(-r,r)$.
		For $0<t<r$, put
		\[
		a(t)=\frac1{p(t)}
		=\frac12+\bigl(\log(1/t)\bigr)^{-\alpha}.
		\]
		The function $a$ is increasing and concave, since
		\[
		a'(t)=
		\frac{\alpha}{
			t\bigl(\log(1/t)\bigr)^{\alpha+1}},
		\]
		and
		\[
		a''(t)=
		-\frac{
			\alpha\bigl(\log(1/t)-\alpha-1\bigr)
		}{
			t^2\bigl(\log(1/t)\bigr)^{\alpha+2}}
		<0.
		\]
		
		Let $I=(x,y)\subset(0,r)$ and $D=y-x$.
		If $D\le x$, concavity yields
		\[
		a(y)-a(x)
		\le D a'(x)
		\le D a'(D)
		=
		\frac{\alpha}{
			\bigl(\log(1/D)\bigr)^{\alpha+1}}.
		\]
		Since $p$ is decreasing on $(0,r)$,
		\eqref{examples-small-interval} implies
		\[
		\frac{
			\|\chi_I\|_{p(\cdot)}
			\|\chi_I\|_{p'(\cdot)}
		}{D}
		\le
		\exp\!\left(
		(a(y)-a(x))\log(1/D)
		\right)
		\le e^\alpha.
		\]
		
		To treat intervals close to the origin, we claim that
		\begin{equation}\label{examples-origin-norm}
			\|\chi_{(0,y)}\|_{p(\cdot)}
			\le 3y^{1/p(y)},
			\qquad 0<y\le r.
		\end{equation}
		Set $k=\log(1/y)$ and
		\[
		A(u)=\frac12+u^{-\alpha},
		\qquad
		\lambda=y^{1/p(y)}=e^{-kA(k)}.
		\]
		The substitution $t=e^{-u}$ gives
		\[
		\int_0^y \lambda^{-p(t)}\,dt
		=
		\int_k^\infty
		\exp\!\left(-u+\frac{kA(k)}{A(u)}\right)\,du.
		\]
		For $u\ge k$,
		\[
		A(k)-A(u)
		\le\alpha k^{-\alpha-1}(u-k),
		\qquad A(u)\ge\frac12.
		\]
		Since $k^{-\alpha}\le1/3$, it follows that
		\begin{align*}
			-u+\frac{kA(k)}{A(u)}
			&=
			-(u-k)+\frac{k(A(k)-A(u))}{A(u)}\\
			&\le
			-\bigl(1-2\alpha k^{-\alpha}\bigr)(u-k)\\
			&\le-\frac13(u-k).
		\end{align*}
		Consequently,
		\[
		\int_0^y \lambda^{-p(t)}\,dt\le3.
		\]
		Since $p\ge1$, replacing $\lambda$ by $3\lambda$
		makes this modular at most one.
		This proves \eqref{examples-origin-norm}.
		
		The conjugate exponent $p'$ is increasing on $(0,r)$.
		Therefore, \eqref{examples-characteristic} gives
		\[
		\|\chi_{(0,y)}\|_{p'(\cdot)}
		\le y^{1/p'(y)}.
		\]
		Combining the last two estimates, we obtain
		\begin{equation}\label{examples-origin-product}
			\|\chi_{(0,y)}\|_{p(\cdot)}
			\|\chi_{(0,y)}\|_{p'(\cdot)}
			\le3y.
		\end{equation}
		If $I=(x,y)\subset(0,r)$ and $D>x$, then $y<2D$.
		By monotonicity of the norms and
		\eqref{examples-origin-product}, the product of the
		two norms over $I$ is at most $6D$.
		
		Evenness handles intervals in $(-r,0)$.
		If $I$ crosses the origin, let
		$y=\max\{|\,\inf I\,|,|\,\sup I\,|\}$.
		Then $I\subset(-y,y)$ and $y\le|I|$.
		Evenness and the triangle inequality give
		\[
		\|\chi_I\|_{p(\cdot)}
		\|\chi_I\|_{p'(\cdot)}
		\le
		4\|\chi_{(0,y)}\|_{p(\cdot)}
		\|\chi_{(0,y)}\|_{p'(\cdot)}
		\le12|I|.
		\]
		This proves the required estimate for all intervals
		contained in $(-r,r)$.
		
		We now check that the estimate persists under the
		constant extension.
		The function $1/p$ is Lipschitz on each of the rays
		$[r/2,\infty)$ and $(-\infty,-r/2]$; let $L_\alpha$
		be a common Lipschitz constant.
		If an interval $I$ is not contained in $(-r,r)$ and
		$D=|I|<r/2$, then it lies in one of these rays.
		Hence \eqref{examples-small-interval} gives
		\[
		\frac{
			\|\chi_I\|_{p(\cdot)}
			\|\chi_I\|_{p'(\cdot)}
		}{D}
		\le D^{-L_\alpha D}
		\le e^{L_\alpha/e}.
		\]
		For $r/2\le D\le1$, the bounds $6/5\le p\le2$ yield
		\[
		\frac{
			\|\chi_I\|_{p(\cdot)}
			\|\chi_I\|_{p'(\cdot)}
		}{D}
		\le D^{-1/3}
		\le(r/2)^{-1/3}.
		\]
		
		Finally, if $D\ge1$, then $p\ge6/5$ implies
		\[
		\|\chi_I\|_{p(\cdot)}\le D^{5/6}.
		\]
		Since $p'=6$ outside $(-r,r)$ and $p'\ge2$ everywhere,
		\[
		\int_I
		\left(\frac1{2D^{1/6}}\right)^{p'(x)}\,dx
		\le
		2^{-6}+\frac{2r}{4D^{1/3}}
		\le2^{-6}+\frac r2<1.
		\]
		Thus
		\[
		\|\chi_I\|_{p'(\cdot)}\le2D^{1/6},
		\]
		which completes the proof that $p\in K_0$.
		The inclusion \eqref{examples-K0-inclusion} now
		gives $p\in\Kplus$.
	\end{proof}
	
	\begin{example}[A one-sided exponent outside $K_0$]
		\label{strict-one-sided}
		Define
		\[
		p(x)=
		\begin{cases}
			2,   &x\le0,\\
			3,   &0<x\le1,\\
			4-x, &1<x<2,\\
			2,   &x\ge2.
		\end{cases}
		\]
		Then
		\[
		p(\cdot)\in\Kplus\cap\Ninf,
		\qquad
		p(\cdot)\notin K_0.
		\]
		In particular,
		\[
		K_0\cap\Ninf
		\subsetneq
		\Kplus\cap\Ninf.
		\]
	\end{example}
	
	\begin{proof}
		Clearly $2\le p(x)\le3$.
		Since $p=2$ outside $(0,2)$, the condition $\Ninf$
		holds with $p_\infty=2$.
		
		To prove $\Kplus$, observe that $x\mapsto p(x)+x$
		is nondecreasing. Consequently,
		\begin{equation}\label{examples-one-sided-drop}
			p(x)-p(y)\le y-x,
			\qquad x<y.
		\end{equation}
		Fix $a<b<c$, and write
		\[
		I^-=(a,b),\qquad I^+=(b,c),\qquad D=c-a.
		\]
		
		Suppose first that $D\le1$, and set
		\[
		P=\operatorname*{ess\,sup}_{I^-}p,
		\qquad
		Q=\operatorname*{ess\,inf}_{I^+}p.
		\]
		By \eqref{examples-one-sided-drop}, $P-Q\le D$.
		Using \eqref{examples-characteristic} and
		$|I^\pm|\le D$, we obtain
		\[
		\|\chi_{I^-}\|_{p(\cdot)}
		\le D^{1/P},
		\qquad
		\|\chi_{I^+}\|_{p'(\cdot)}
		\le D^{1-1/Q}.
		\]
		If $P\le Q$, their product is at most $D$.
		If $P>Q$, then $PQ\ge4$ and
		\[
		\frac1P-\frac1Q
		=-\frac{P-Q}{PQ}\ge-\frac D4.
		\]
		Thus, in either case,
		\[
		\frac{
			\|\chi_{I^-}\|_{p(\cdot)}
			\|\chi_{I^+}\|_{p'(\cdot)}
		}{D}
		\le D^{-D/4}
		\le e^{1/(4e)}.
		\]
		
		Suppose next that $D\ge1$.
		Since $p(\cdot)\ge2$,
		\[
		\int_{I^-}
		\left(\frac1{\sqrt D}\right)^{p(x)}\,dx
		\le\frac{|I^-|}{D}\le1,
		\]
		and hence
		\[
		\|\chi_{I^-}\|_{p(\cdot)}\le\sqrt D.
		\]
		Moreover, $p'(\cdot)=2$ outside $(0,2)$ and $p'(\cdot)\ge3/2$
		everywhere. Therefore,
		\[
		\int_{I^+}
		\left(\frac1{4\sqrt D}\right)^{p'(x)}\,dx
		\le
		\frac{|I^+\setminus(0,2)|}{16D}
		+\frac{|I^+\cap(0,2)|}{(4\sqrt D)^{3/2}}
		\le\frac1{16}+\frac14<1.
		\]
		It follows that
		\[
		\|\chi_{I^+}\|_{p'(\cdot)}\le4\sqrt D.
		\]
		Combining the two cases proves $p\in\Kplus$.
		
		To see that $p\notin K_0$, let $0<h<1$ and
		$I_h=(-h,h)$.
		The exponent equals $3$ on $(0,h)$, whereas its
		conjugate equals $2$ on $(-h,0)$.
		By monotonicity of the Luxemburg norm,
		\[
		\|\chi_{I_h}\|_{p(\cdot)}
		\ge\|\chi_{(0,h)}\|_{p(\cdot)}=h^{1/3},
		\]
		and
		\[
		\|\chi_{I_h}\|_{p'(\cdot)}
		\ge\|\chi_{(-h,0)}\|_{p'(\cdot)}=h^{1/2}.
		\]
		Consequently,
		\[
		\frac{
			\|\chi_{I_h}\|_{p(\cdot)}
			\|\chi_{I_h}\|_{p'(\cdot)}
		}{|I_h|}
		\ge\frac12h^{-1/6}
		\longrightarrow\infty
		\qquad\text{as }h\downarrow0.
		\]
		Thus $p(\cdot)\notin K_0$.
	\end{proof}
	
	\begin{remark}
		By Theorem~\ref{main-Theorem}, the right-sided maximal
		operator $M^+$ is bounded on $L^{p(\cdot)}(\R)$
		for the exponents in
		Examples~\ref{non-log-holder}
		and~\ref{strict-one-sided}.
		The last example shows that the one-sided criterion
		applies to exponents excluded by the classical
		$K_0$ condition, even under the additional
		assumption $\Ninf$.
	\end{remark}
	
	\section*{AI Disclosure}
	The original ideas underlying this work are the author’s own. During the preparation of this manuscript, the author used Claude and ChatGPT to discuss technical details of the proofs, assist in drafting and checking parts of the mathematical arguments, organize references, and improve the clarity and readability of the exposition. All mathematical statements and proofs in the final manuscript were written, reviewed, and verified by the author, who takes full responsibility for the content of this manuscript.
	

\end{document}